\documentclass[12pt,amscd,amssymb]{article}
\usepackage[dvips]{graphics}
\title{
Numbers as Functions
(The Development of an Idea in the Moscow
School of Algebraic Geometry)\footnote{Published in: Bolibruch, A.A. (ed.) et al., Mathematical events of the twentieth century. Berlin: Springer; Moscow: PHASIS, 2006, 297-329. I am very much grateful to Yu. N. Torkhov for  preparation of the TeX files for  pictures.}}
\author{A.N.Parshin}
\date{}

\def\F1{{\bf F}_1}

\def\Spec{\mbox{Spec}}
\def\Gal{\mbox{Gal}}
\def\exp{\mbox{exp}}
\def\log{\mbox{log}}
\def\res{\mbox{res}}
\def\deg{\mbox{deg}}
\def\Prod{\prod\limits}

\def\P{{\bf P}^1}
\def\R{{\bf R}}
\def\R{{\bf R}}
\def\Z{{\bf Z}}
\def\Q{{\bf Q}}

\def\C{{\bf C}}
\def\Fp{{{\bf F}_p}}
\def\Fq{{{\bf F}_q}}

\def\Qp{{{\bf Q}_p}}

\begin{document}
\maketitle

Where numbers come from nobody knows. Ethnographers have traveled through all the countries of the world, up and down, backward and forward, and have found people for whom "one", "two", and "many" are sufficient. And yet, these people have refined arts, subtle mythology, and nontrivial crafts. They are people quite as much as we are, but without that "one, two, three", and so on to infinity. Prometheus did not bring them "the science of number, the most important of all sciences".

However, there is no need to travel the world to see the sharp difference between the first numbers and those that followed. Language has retained enough evidence of that. Thus, in many linguistic families the etymology and grammatical forms of the first three or four cardinal numbers are fundamentally different from those of all other numbers. Moreover, among all the peoples of the world the "first" numbers are burdened with a rich symbolism and have their own individual character. In the sterile series of natural numbers all this disappears completely.

That there are infinitely many numbers seems to have been recognized for the first time in ancient Greece. EuclidТs {\em Elements} even contain a proof that
the series of prime numbers is infinite. Here the infinite is understood as a potential infinity, a non-finiteness. To the modern person the origin of numbers is completely comprehensible: they arose from counting "things' or "objects" (but where does counting come from?). It also seems obvious to modern people that, once having begun to count, they are unable to stop. To imagine a finite closed universe of numbers in real life is, after all, not easy. Still, in mathematics there are, for example, finite fields.

Whatever the situation with numbers, they are originally a discrete, countable object. All those irrationalities and continua that the Greeks struggled with arose later in history. But the concept of a function seems to have arisen as the incarnation of something continuous, the trajectory of a stone that is thrown, or a line drawn in the sand with a finger. Functions are connected with motion.

However, the subsequent evolution of algebra and many functions turned into something discrete, amenable to algorithmization by means of some sort of Maple V or PARI. Many have speculated on the relation of the discrete and the continuous in mathematics. Hermann Weyl wrote about these two modes of understanding. And Andre Weil related the following УargumentФ Claude Chevalley and Oscar Zariski. What is a curve? They went to the board and wrote the following:
\par\smallskip
\includegraphics{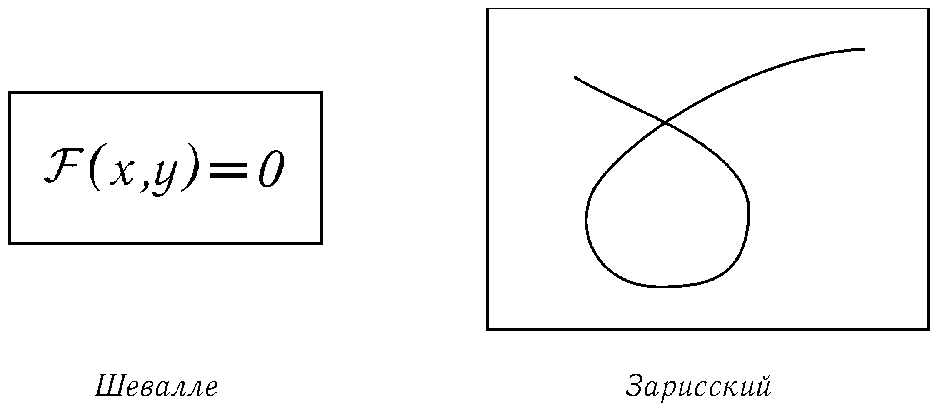}
%$$
%\begin{array}{ll}
%\mbox{Chevalley} &   f(x, y) = 0 \\
%\mbox{Zariski} &   \mbox{Ёшёєэюъ}
%\end{array}
%$$
\par\smallskip
We can see both answers. Both are generated by a sweep of the hand, but the formula can also be read aloud. Thus we are faced with the ancient quarrel between the ear and the eye, the world of language and the world of vision\footnote{The current generation adds here the computer keyboard and the mouse.}.

The analogy between numbers and functions that forms the subject of this article is an even greater leap . . . In this article the continuous and the discrete enter on both sides Ч into numbers and into functions. Algebra and analysis work hand in hand here.

The basis of our discussion is a lecture delivered by the author at the conference "Mat\'eriaux pour
l'Histoire des Math\'ematiques au XX\`eme si\`ecle",  which
took place at LТUniversite de Nice Sophia-Antipolis in Nice in January 1996. My task was to describe one area in the development of arithmetical algebraic geometry in Moscow during the 1950s and 1960s. I made no attempt to present
any complete historical study of the development of algebraic geometry during this "golden age of the Moscow mathematical school."

We shall begin by explaining the meaning of the analogy between numbers and functions, starting with the simplest concepts. In the second part we study a nontrivial example: the explicit formula for the law of reciprocity. In the third part we shall become acquainted with certain aspects of the "social" life of the Moscow school, in particular, with certain seminars, lectures, and books. In the final part we shall examine another example of this analogy: arithmetical surfaces, an example that is indisputably the summit of this area. As for the time frame, I shall hardly pass beyond the early 1970s.

The fact is that the smooth development of this idea Ч the analogy between numbers and functions Ч which began in the last third of the nineteenth century, underwent a sharp jump in the 1960s. It was recognized that the preceding development had occurred in the framework of a one-dimensional world. It became clear that one could and should pass to other dimensions. How exactly this jump took place we wish to relate here.

The interested reader may consult\cite{ALNT},
  \cite{ZP}, \cite{S1,S2} to get acquainted with the subsequent development of these ideas, which tended to be broad rather than deep. These same articles contain some results that we have omitted.

  \section{The Analogy}

To understand the origin of the analogy between numbers and functions, let us look at the following table:
     $$
     \begin{array}{ccc}
     f \in \Z \subset \Q & F \in \Fp[t] \subset \Fp(t)  &
F \in \C[t]  \subset \C(t)\\
          && \\
f = (\pm)p_1^{\nu_1} \dots p_n^{\nu_n}  &

F  =  aP_1^{\nu_1}  \dots P_n^{\nu_n}  &

F  =  a(t-t_1)^{\nu_1} \dots (t-t_n)^{\nu_n}   \\
               && \\
f \ne 0 &  F \ne 0,~a \in \Fp^{*}  &  F \ne 0,~a \in \C^{*}
\end{array} $$

Here we are comparing the ring $\Z$ of integers and the rings of polynomials $\Fp[t]$ and $\C[t]$ in one variable $t$ (with coefficients, of course, in the finite field $\Fp$ consisting of $p$ elements and the field $\C$ of complex numbers)\footnote{Here and below $k^{*}$ is the set of nonzero elements of the field $k$, that is, its multiplicative group.}. The nonzero elements of these rings (the numbers $f$ and functions $F$) can be expanded as the product of prime numbers $p$ and irreducible polynomials respectively. Over the field $\Fp$ the latter correspond to conjugate elements of the algebraic closure of the field $\Fp$. Over the field $\C$ they are linear polynomials $t-t_0$, where уфх $t_0$ is an arbitrary point on the complex line $\C$.

The integers $\nu_k$ that occur in the given expansion have the following fundamental property: they are valuations \footnote{Such valuations are usually called non-archimedean. Corresponding to them are the multiplicative (non-archimedean) norms $\vert f \vert$, i.e. $\vert f \vert = p^{-\nu (f)}$, for which $\vert fg \vert = \vert f \vert
g \vert$  and $\vert f+g \vert \le max (\vert f \vert , \vert g \vert)$.  If this last condition is weakened to $\vert f+g \vert \le \vert f \vert + \vert g \vert$, we obtain the well-known definition of a norm.}
in the sense that
\begin{itemize}
     \item $\nu(fg) = \nu(f) + \nu(g)$
     \item $\nu(f + g) \ge min(\nu(f), \nu(g))$.
     \end{itemize}

For our rings the valuations $\nu$ assume nonnegative values, and they can be uniquely extended to the quotient-fields (the field $\Q$ of rational numbers and the fields $\Fp(t)$ and $\C(t)$ of rational functions) as homomorphisms onto the group $\Z$. Accordingly, in these fields we have expansions of their elements into products that generalize the expansions considered above.

This is the first observation showing that number rings and rings of functions have certain properties in common.
We now call attention to the fact that in the case of $\C$ the set of valuations coincides completely with the set of points of the complex affine line. The same is true for the affine line over the finite field $\Fp$ if such points are taken as the maximal ideals of the ring $\Fp[t]$. Each such ideal is a principal ideal, that is, consists of the multiples of some irreducible polynomial $P$. In this case, the base field is not algebraically closed, and the correct definition of the points differs from the straightforward one: a point is characterized by its coordinate, that is, an element of the algebraic closure of the field $\Fp$.

We may also attempt to use our geometric intuition and introduce a geometric object in the case of the ring $\Z$. We denote it by $\Spec(\Z)$, and we shall at first regard it as the set of prime numbers $p = 2,3,5,\dots$ or prime ideals $(p) \subset \Z$. In this way our table expands to the following:
$$
     \begin{array}{ccc}
     \nu_p  &  \nu_P    &   \nu_P \\
            && \\
     p \in \Spec(\Z) & P \in \Spec(\Fp[t])  & P = (t=t_0) \in \Spec\C[t]) \\
                 && \\
                   &  \mbox{affine line over}~\Fp  &  \mbox{affine line over}~\C
     \end{array}
     $$

For any point $P$ on the affine line over $\C$ we have the power-series expansion of a rational function $F$:
$$
     F = \sum_{i_0}^{\infty} a_i(t-t_0)^i, \mbox{уфх}~i_0 =
     \nu_P(F), ~a_i \in \C.
     $$

There is also an expansion of this type in the case of the line over $\Fp$. The analogous construction for the field $\Q$ is the $p$-adic representation of rational numbers:
$$
     f = \sum_{i_0}^{\infty}  a_ip^i, \mbox{уфх}~i_0 = \nu_p(f)~
     \mbox{ш}~a_i \in  {0, 1, \dots, p-1}.
     $$
Both expansions are connected with field embeddings: $\Q$ into the field $\Qp$ of $p$-adic numbers, and $\Fp(t)$ and $\C(t)$ into the fields of power series $\Fp((t))$ and $\C((t))$. These embeddings are the completions of the fields relative to the metrics defined by the valuations(see \cite{BS, CF}):
$$
  \rho(x,y) = p^{-\nu(x - y)}
$$
in the case of the fields $\Qp$ or $\Fp((t))$, and
$$
  \rho(x,y) = c^{-\nu(x - y)}
$$
in the case of $\C((t))$. (Here $c \ne 0$ is an arbitrary fixed constant.)

As numerical variants of power series the $p$-adic numbers were introduced by Kurt Hensel \cite{Hen}. The analogy between power series and the expansions of rational numbers in powers of $p$~ (for $p = 10$) had been considered earlier by Isaac Newton.
A more profound manifestation of this analogy is the global property of valuations known in number theory as {\em the product formula}. To obtain it we must enlarge our objects in order to make them "compact" or "complete". In the case of the affine line it is necessary to embed it in the projective line $\P1$ by adding a point at "infinity". It corresponds to the valuation
$$
     \nu_{\infty}(f) = \deg(f).
    $$
The point at "infinity" has no meaning as an ideal of the ring of polynomials in $t$. (The set of such ideals is exhausted by the points of the affine line.) However, our projective line $\P1$ contains another affine line(the complement of the point $t = 0$), which corresponds to the subring $\Fp[t^{-1}]$ or $\C[t^{-1}]$ of the field of rational functions. And the "infinite" point corresponds to the ideal $(t^{-1})$ of this ring. Thus, in the case of a field of functions all points are arranged Уin the same wayФ: they correspond to nonarchimedean valuations of the field of rational functions, and in this way we obtain all valuations of the field.

In the case of the field $\Q$ our geometric object $\Spec(\Z)$ is also not "compact". The prime numbers correspond to all nonarchimedean valuations of the
field $\Q$.
But there is also an archimedean valuation
\footnote{
If $K$ is any field, we shall take the (archimedean) valuation to be $\log\vert f \vert$, where $\vert f \vert$ is a norm on $K$.}
     $$
     \nu_{\infty}(f) = - \log \vert f \vert,~f \in \Q^{*}.
     $$

and by a theorem of Ostrovskii we now have all valuations of the field $\Q$. The fundamental difference with the geometric case is that in the number field situation the point at "infinity" has no meaning as an ideal of some subring of the field $\Q$.
The product formula for the field $\Q$ has the form
$$
     (\prod_{p \in \Spec(\Z)}p^{-\nu(f)})\times \vert f \vert  = 1,~f
     \in \Q^{*}.
     $$

To compare it with the corresponding formula for the projective line, we pass from the product to the sum
\begin{displaymath}
     \sum_{p \in \Spec(\Z)} \nu_{p}(f)\mbox{log}p +  \nu_{\infty}(f)  = 0.
     \end{displaymath}
 For the projective line $X$ over $\Fp$ we have
\begin{displaymath}
     \sum_{P \in X} \nu_{P}(f)\mbox{deg}P +  \nu_{\infty}(f)  = 0,
     \end{displaymath}

and for the projective line over $\C$
\begin{displaymath}
     \sum_{P} \nu_{P}(f) +  \nu_{\infty}(f)  = 0.
     \end{displaymath}
     This means that the polynomial $f$ has a number of zeros equal to its degree.

The projective line is a special case of an algebraic curve, and the ring $\Z$ is a special case of rings of integers in fields of algebraic numbers(finite extensions of the field $\Q$). These two concepts combine in the language of the theory of schemes as schemes of dimension 1.

A scheme is a space with a sheaf of rings: the structure sheaf ${\cal O}_X$ of regular functions on the scheme $X$. For each open set $U \subset X$ we know the set of functions that are regular on $U$ (namely ${\cal O}_X(U)$). In algebraic geometry the fiber of the sheaf ${\cal O}_{X,x}$ at the point $x \in X$ consists of the rational functions that do not have a pole at that point. In the example $X = \Spec(\Z)$ that we have considered one may set
$$
{\cal O}_{X,p} = \{f \in \Q : f = m/n~\mbox{ё}~m, n \in \Z~\mbox{ш}
    ~(n,p) = 1 \}
$$

In general schemes $X$ (of finite type) over the ring $\Spec(\Z)$ are the basic object of arithmetical algebraic geometry. There are two types of schemes. Roughly speaking, they are УsetsФ defined by equations with finite coefficients and sets defined by equations with integer coefficients. We shall denote these two cases below as {\em geometric} (or {\em functional}) and {\em arithmetical} respectively.

Our basic examples $\Spec(\Fp(t))$ and $\Spec(\Z)$ happen to be the simplest representatives of these two types. The original classification of schemes consists of classification according to dimension. By that we mean the absolute dimension over the ring $\Z$. For affine schemes it coincides with the Krull dimension of the corresponding ring (that is, the length of a maximal chain of prime ideals\footnote{We recall that an ideal $\wp$ of a ring $A$ is prime if the quotient ring $A/\wp$ has no zero divisors and is not the zero ring {0}. That is, the ring $A$ itself is not a prime ideal.}).

The examples with which we began our exposition happen to be schemes of dimension 1. Ordered chains of prime ideals in $\Z$, $\Fp[t]$ and $\C[t]$ have length 1. For example, in $\Z$ we have $(0) \subset (p)$, and in $\C[t]$ we have $(0) \subset (t-t_0)$.
From the point of view of arithmetic, the ring $\C[t]$ is not of arithmetic type. We, however, have included it in our picture as the example of a geometric object closest to our intuition.
One of the routes in arithmetic consists of a transition from varieties over the field $\C$ to varieties over $\Fq$ and then to schemes over $\Spec(\Z)$. Such an approach suggests the correct statements of theorems that are valid in both situations and sometimes also methods of proving them.

The terminology of schemes arose only in the mid-twentieth century, but attempts to unify number theory and algebraic geometry into a single subject had been made much earlier. Probably, the first person to recognize the importance of the concept of dimension for arithmetic was Kronecker. As early as in the nineteenth century he attempted to develop arithmetic not only for dimension 1 but also for arbitrary dimensions. This program was neglected for a long time and was resurrected only in the mid-twentieth century.
We can point to two ground-breaking lectures at International Congresses of Mathematicians in which this problem was discussed. The first was by
A.Weil at the 1950 Congress in Cambridge, Massachusetts \cite{W1}. Weil described Kronecker's goals as follows: "He was, in fact, attempting to describe and to initiate a new branch of mathematics, which would contain both number theory and algebraic geometry as special cases".
The second lecture was by I.R.Shafarevich and was given at the Stockholm Congress in 1962 \cite{Sh2}.

Between these two events A.Grothendieck created the theory of schemes (\cite{Gr, D3}). I think that Weil's lecture had some influence on Grothendieck. As for ShafarevichТs lecture, he was now able to use the language of schemes as the foundation for the further development of arithmetic.
Using the concept of a scheme we can describe our analogy by the following table, where we compare schemes $X$ of the same dimension from both parts of the table:
\vskip 0.5cm

\begin{tabular}{|c|c|c|}
\hline
 $\mbox{dim}(X)$ & geometric case & arithmetical case \\
           ...     &          ...            &                 ...        \\
&& \\
\hline
2  & algebraic surfaces/$\Fq$ &  arithmetical surfaces \\
&& \\
\hline
1    &  algebraic curves/$\Fq$       &  arithmetical curves  = \\
  &                         &    finite coverings of  $\Spec(\Z)$ \\
&& \\
\hline
 0    &        $\Spec(\Fq)$     &   $\Spec(\F1)$ \\
\hline
\end{tabular}
\vskip 0.5cm
Here $\F1$ is the УfieldФ of one element (see below).

This table is the result (or starting point) of a completely new way of looking at the analogy between numbers and functions. Over a period of almost 80 years only the row of the table relating to dimension 1 was known and studied.

The leading role in this development belonged to D.Hilbert \cite{H1, H3, W}. This analogy was one of his favorite ideas, and it was thanks to Hilbert that it achieved fame and became one of the central ideas in the development of number theory during the twentieth century.

In this section and those that follow we shall speak of this Hilbert period, then pass to the description of the jump to other dimensions that occurred in the 1960s.
Let us now repeat the constructions given above in the more general situation of arbitrary curves (or schemes) of dimension 1.

Let $X$ be an algebraic curve over a finite field $\Fq$, let $K = \Fq(X)$ be the field of rational functions on $X$, and let $\nu_{x}: K^{\star}\rightarrow \Z $ be the valuations corresponding to closed points $x \in X$. If we assume that $X$ is a {\em projective} curve, we have the "sum formula"
\begin{displaymath}
\sum_{x \in X} \nu_{x}(f)\mbox{log}\# k(x) = 0,~f \in K^{\star},
\end{displaymath}
or the product formula
\begin{displaymath}
\prod_{x \in X} \vert f \vert_{x} = 1,
\end{displaymath}
where
\begin{displaymath}
\vert f \vert_{x} = \# k(x)^{-\nu_{x}(f)} \quad ,
\end{displaymath}
Here $k(x) = {\cal O}_{X,x}/m_x$ is the field of residues of the local ring ${\cal O}_{X,x}$ at the point $x \in X$, $m_x$ is a maximal ideal, and $k(x)$ is a finite extension of the field $\Fq$. In the geometric case we can use either the curve $X$ itself (the point of view of algebraic geometry) or the field $K$ of rational functions on X (the point of view of algebra). According to a well known result, these are two descriptions of the same object. (Every field of algebraic functions of one variable has a unique projective nonsingular curve X as a model.)

If we turn to arithmetic, we can observe that the algebraic point of view was dominant for a long time. The object of study was a field $K$ of algebraic numbers, that is, a finite extension of the field $\Q$ of rational numbers. But now we can also use the geometric point of view, that is, the point of view of the theory of schemes. This has the following appearance.

Let $X$ be the set of prime ideals $\wp$ of the ring $\Lambda_K$ of integers in the field $K$. To each $\wp \in X$ there corresponds a valuation $\nu_{\wp}$, namely
$$
\nu_{\wp}(f) = log \vert f \vert_{\wp},~f \in \Q^{*},
$$
where
$$
\vert f \vert_{\wp} = \# (\Lambda_K /\wp)
$$
is the corresponding norm. It is easy to see that for $K = \Q$ this definition coincides with the previous definition.

As before, the product of $\vert f \vert_{\wp}$ over all $\wp$ is again not equal to 1. But now we must adjoin a finite number of УinfiniteФ points $\infty$, where $\infty$ is a certain embedding of the field $K$ in the field $\C$ of complex numbers. The number of such embeddings equals the degree $[K:\Q]$  of the field $K$ over $\Q$. If the embedding $\infty$ is real, that is, has the form $K \subset \R \subset
\C$, the norm equals
$$
\vert f \vert_{\infty} = \vert f \vert_{\R}.
$$
Otherwise, we have
$$
\vert f \vert_{\infty} = \vert f \vert_{\C} = \vert f \vert^2.
$$
We then have the product formula:
\begin{displaymath}
\prod_{x \in X' \cup \infty} \vert f \vert_{x} = 1 \quad .
\end{displaymath}

Of course, all these valuations have a simple interpretation. They correspond to {\em all} possible completions of the field $K$, namely, $\wp$--adic fields $K_{\wp}$, real fields $\R$  and complex fields $\C$. For the field $\Q$ we have a unique embedding $\Q \subset \R$.

A scheme structure on $X$ is defined by the sheaf ${\cal O}_X$ whose fibers are
$$
     {\cal O}_{X,\wp}
= \{ f \in K : \nu_{\wp}(f) \ge 0 \}.
     $$
The rings ${\cal O}_{X,\wp}$ contain a maximal ideal  $m_x =
\{ f \in K : \nu_{\wp}(f) \ge 0 \}$, and completing with respect to it, we obtain a complete local ring $\hat{{\cal O}}_{X,\wp}$. Its field of fractions will be the completion of $K$ with respect to the valuation $\nu_{\wp}$. For УinfiniteФ points there is no such construction, there are only the fields $\R$ and $\C$, but no subrings in them. For that same reason we cannot introduce the structure of a scheme on the entire set $X \cup \{\infty\}$ of УpointsФ of the field $K$.

A field extension $K \supset \Q$ gives a mapping of degree $[K : \Q]$
$$
  X \cup \{\infty\} \rightarrow \Spec (\Z) \cup \infty
$$
and above an infinite point of the scheme $\Spec(\Z) \cup \infty$ lie exactly $[K :\Q]$ infinite points of the scheme $X \cup \{\infty\}$.

But we want to move in the opposite direction, from geometry to arithmetic. And the theory of schemes makes it possible for us to apply the language of geometry in the situation of number theory. If $X = \Spec(\Z)$, the closed points $x$ of the scheme $X$ are the primes $p$, and we have a canonical isomorphism:
\begin{displaymath}
k(x) \cong \Fp \quad .
\end{displaymath}
Here $ k(x) = {\cal O}_{x}/m_{x}$, where $m_{x}$ is the maximal ideal of
${\cal O}_{x}$.

We can speak of rational {\em numbers} $f\in \Q$ as rational {\em functions} on the "curve" $X$ with values $f(x)\in k(x)$. The fundamental difference from genuine curves is that the values $f(x)$ of our function belong to different fields $\Fq$, as $x$ ranges over the УcurveФ $X$. The fields $\Fq$ are different from each other. They are not extensions of the same finite field $\Fp$, as was the case with curves. They contain as a common subfield only the УfieldФ $\F1$ consisting of one element. We include it in our table under dimension 0 as the final object in the category of schemes of arithmetical type\footnote{Surprisingly, this is not a vacuous concept. It has a rich structure. For example, the higher $K$-groups $K_{.}({\bf F}_{1})$ are defined, and they coincide with the stable homotopy groups of spheres (See \cite{S4}).}

\section{The Reciprocity Law}

Up to now we have spoken only about the simplest aspects of the analogy between numbers and functions. A much more profound fact is the product formula for the normed residue symbol
$$ \left( \frac{\lambda, \mu}{\wp} \right),$$
discovered by Hilbert. In \cite{H1} he wrote, "The reciprocity law
$$ \Prod_{\wp} \left( \frac{\lambda, \mu}{\wp} \right) = 1$$
reminds  the Cauchy integral theorem, according to which the integral of a function over a path enclosing all of its singularities always yields the value 0. One of the known proofs of the ordinary quadratic reciprocity law suggests an intrinsic connection between this number-theoretic law and CauchyТs fundamental function-theoretic theorem."\footnote{"Das Reziprozit\"atsgesetz in der Fassung
$$ \Prod_{\wp} \left( \frac{\lambda, \mu}{\wp} \right) = 1$$
 erinnert an den Cauchyschen Integralsatz in der Funktionentheorie, dem zu folge ein complexes Integral,um alle einzelnen Singularit\"aten einer Funktion gef\"uhrt, insgesamt stets den Wert 0 gibt.
 Einer der bekannten Beweise des gewohnlichen quadratischen Reziprozit\"atsgesetzes weist auf
einen inneren Zusammenhang zwischen jenem zahlentheoretischen Gesetz und Cauchys funktionentheoretischen Fundamentalsatz hin." (David Hilbert. {\em Gesammelte Abhandlungen}, Erster Band, {\em Zahlentheorie}, NewYork:Chelsea,1965(reprint),pp.367-368.) - Transl.}

This idea was realized by Shafarevich in his purely local construction of the Hilbert symbol. The proof of reciprocity given by him was a far-reaching extension of the corresponding result for residues \cite{Sh1}. This result is an important contribution to the solution of Hilbert's ninth problem. (See the statement of it in \cite{H4, Fad}; and see commentaries on it in the Russian edition of \cite{H3}. Shafarevich was probably the first in our country to take this analogy seriously.

He used it in a highly non-trivial manner, since it was necessary to compare $p$-adic number fields whose multiplicative groups had a complex structure with the much simpler fields of power series. ShafarevichТs paper begins with the quotation from Hilbert given above. He then corrects Hilbert, showing that the analog of the product formula must be a formula for the sum of the residues rather than the Cauchy integral theorem.

We first recall the well known constructions from class-field theory. Class-field theory is a method of describing Abelian extensions of a field $K$ of arithmetic type, such as $\Q$ or $\Fp(t)$ (that is, extensions $L/K$ with a commutative Galois group $\Gal(L/K))$. In this case it is called global class-field theory.
If $K$ is a number field, it can be embedded in the completion $K_{\wp}$ for all prime ideals $\wp$, as we saw above. In this section we shall be dealing only with the fields $K_{\wp}$. The field $K_{\wp}$ is called a local field, and describing its Abelian extensions is a local class-field theory problem. To this end, let us consider a maximal Abelian extension $K_{\wp}^{ab}$
 as the union of all finite Abelian extensions. The problem is to describe its Galois group over $K_{\wp}$ using a construction intrinsically connected with the field $K_{\wp}$ rather than with its extensions.

The main result of local class-field theory is the existence of a canonical homomorphism
$$
     \varphi : K_{\wp}^{*} \rightarrow \Gal(K_{\wp}^{ab}/K_{\wp}),
     $$
which has trivial kernel and dense image. Global class-field theory for the field $K$ then reduces in a natural way to the local theories for all the fields $K_{\wp}$ (see, for example, \cite{CF}).
We shall show how the mapping $\varphi$ defines the Hilbert symbol. Let us assume that a  $p^n$-th root of unity $\zeta$  belongs to our field. Here $p$ is the characteristic of the field of residues. We take two numbers $\lambda$ and $\mu$ from $K_{\wp}$.  In this situation
we have an Abelian extension $K(x)/K$, where $x^{p^n} = \lambda$. Its Galois group $G$ is the cyclic group of order $p^n$, and for every
$\sigma \in G$  we have
     $$
     \sigma(x) = (\mbox{some power of the root}~\zeta)x.
$$
From class-field theory we obtain a mapping
$$
     K_{\wp}^{*} \rightarrow \Gal(K_{\wp}^{ab}/K_{\wp}) \rightarrow
      G,
     $$

which we also denote $\varphi$. We can now define the Hilbert symbol by the condition
$$
     (\sqrt[p^n]{\lambda})^{\varphi(\mu)}  = \left( \frac{\lambda, \mu}{\wp}
     \right) \sqrt[p^n]{\lambda},
     $$

where $\sqrt[p^n]{\lambda} = x$, and the result is independent of the choice of $x$. Thus, to define this symbol, we must {\em leave} our local field and work with its extensions. The problem posed by Hilbert was to obtain an explicit expression entirely within the field $K_{\wp}$ and then use it to reverse this process by constructing the mapping $\varphi$ and developing class-field theory.

Further, if we take $\lambda$ and $\mu$ from the original global field $K$ rather than the local field, we obtain symbols for all prime ideals $\wp$. To obtain a global reciprocity law, one must also define a symbol for infinite points $\infty$. If we also define a symbol for them (which is much simpler to do), we obtain the reciprocity law described by Hilbert (see above).

In particular, let $K = \Q, ~p= 2,~n = 1$, and take as $\lambda, \mu$ two odd primes $a$ and $b$. The only factors that remain in the infinite product of the general reciprocity law are those corresponding to $\wp = (a), (b)$ and $\infty$. HilbertТs law then reduces precisely to the quadratic reciprocity law of Gauss
$$
     \left( \frac{a}{b} \right) \left( \frac{b}{a} \right) =
     (-1)^{ \frac{a-1}{2} \cdot \frac{b-1}{2}},
     $$
     where $\left (\frac{a}{b} \right )$ is the Legendre symbol. I now pass to the explanation of ShafarevichТs construction and the way in which it is connected with the residues of differential forms on Riemann surfaces. Shafarevich considered the case $n = 1$. The general case, just like the application to the construction of class-field theory starting with the local definition of the symbols, was considered by A. I. Lapin \cite{La1, La2, La3}\footnote{His first paper was written in 1950, when he was in detention. The question of the possibility of publishing it was discussed in the Central Committee of the Communist Party of the USSR on the request of Academy of Sciences and the Ministry of Internal Affairs. After the question had been considered by three members of the Politbyuro, one of whom was L.P.Beria, permission was given to publish it under a pseudonym. However, by then Lapin had been freed, and the paper was published in the usual way. The materials of this correspondence were recently discovered in the archives of the Central Committee (See Voprosy Istorii Estestvoznaniya i Tekhniki, 2001, No.2, 116-128). This issue also contains reminiscences of S. S. Demidov, I. R. Shafarevich and
     I. G. Bashmakova on A. I. Lapin.}.

For brevity, we shall denote our local field by $K$. It is the field of fractions of a discrete valuation ring ${\cal O}$ with maximal ideal $\wp$ and with the field of residues ${\cal O}/\wp = \Fq$. We denote the generator of the ideal $\wp$ by $\pi$.
The multiplicative group $K^{*}$ has the following structure:
$$
     K^{*} = \{ \pi^m \} {\cal O}^{*} = \{ \pi^m \} RU,
     $$

where the set $R$ consists of multiplicative representatives of the field of residues $\Fq$, and $U = 1 + \wp$ is called the group of principal units.
The Hilbert symbol has two important properties that are useful for its computation, namely
$$
      \left( \frac{\lambda\cdot\lambda', \mu}{\wp} \right) =
      \left( \frac{\lambda, \mu}{\wp} \right)
      \left( \frac{\lambda', \mu}{\wp} \right),
     $$
   and
     $$
      \left( \frac{(\lambda)^{p^n}, \mu}{\wp} \right)  =
      \left( \frac{\lambda, \mu}{\wp} \right)^{p^n}  =  1.
     $$
 The same is true for the second argument $\mu$.

These properties show that to compute this symbol we must find some system of generators for the group $U/U^{p^n}$. (For the group $R$ we have $R = R^{p^n}$). For this purpose Shafarevich used the Artin-Hasse exponentials $E(\alpha,  x)$ and the variant of them $E(\alpha)$. They are defined for elements $\alpha$ of a maximal unramified subring ${\cal O}_{nr} \subset
{\cal O}$ and $x \in \wp$. These functions are homomorphisms from the ring ${\cal O}_{nr}$ into the group of units $U$. We shall see that they are the analogs of the exponential functions. We shall find the following abbreviation useful:
$$
     \lambda \approx \mu \Longleftrightarrow ~\lambda\mu^{-1}
     \in K^{p^n}.
     $$

We have the following fundamental expansions:
 $$
     \lambda \approx \pi^aE(\alpha)\Prod_{\begin{array}{c}
     1 \le i < pe/p-1 \\
     p \not | i
     \end{array}}E(\alpha_i,\pi^i),
     $$
     $$
     \mu \approx \pi^bE(\beta)\Prod_{\begin{array}{c}
     1 \le i < pe/p-1 \\
     p \not | i
     \end{array}}E(\beta_i,\pi^i).
     $$
The integers $a$ and $b$ are defined modulo $p^n$, and the values of all $E$-functions
are defined modulo $p^n$-powers. If we introduce the homomorphism $\delta : U/U^{p^n}
      \rightarrow U/U^{p^n}$ with
     $\delta(\lambda) = E(\alpha), \delta(\mu) = E(\beta)$, the required local expression will have the following appearance:
$$
      \left( \frac{\lambda, \mu}{\wp} \right) =
     E(\beta)^{a}E(\alpha)^{-b}E(\gamma),
     $$
where
$$
E(\gamma) \approx \delta(\Prod_{i,j}E(i\alpha_i\beta_j, \pi^{i+j})).
$$
%(here $\alpha = \sum_{i \ge 0}\alpha_ip^i$, with $\alpha_i \in R$)

The main thing is to show that the result is independent of the choice of the prime element $\pi$. This is true, but the proof is complicated and rather long.

Even so, it remained unclear how to find the value of $\gamma$ explicitly. This was done later by two mathematicians independently of each other, H.Br\"uckner in Germany and S.V.Vostokov in Leningrad \cite{Br1, Vos, Br2}.

To understand the analogy with Riemann surfaces, let us consider a point $P$ on such a surface, a local coordinate $t$, and the corresponding field $K = \C((t))$ of Laurent power series. For the multiplicative group of the field $K$ we have:
$$
     K^{*} = \{t^m\}\C^{*}U
     $$
and all elements $\lambda,\mu \in U$ have an expansion:
     $$
     \lambda = exp(A) = \Prod_{i \ge 1}exp(\alpha_it^i),~\alpha_i \in \C,
     $$
     $$
     \mu = exp(B) = \Prod_{j \ge 1}exp(\beta_jt^j),~\beta_j \in \C.
     $$

In the field $K$ there are two simple operations: taking the derivative $\partial = d/dt$ and the residue $\res(\sum \alpha_it^i) = \alpha_{-1}$.
The analogy with the residue of a differential form at the point $P$ can now be seen from the following table:
$$
     \begin{array}{ccc}
     \exp(\alpha_it^i)  &  \sim  &  E(\alpha_i, \pi^i) \\
&& \\
     \exp(\beta_jt^j)  &  \sim  &  E(\beta_j, \pi^j) \\
     && \\
     exp(B\partial A) & \sim & \Prod_{i,j}E(i\alpha_i\beta_j,
\pi^{i+j})) \\
&& \\
     \res(\exp(B\partial A)) & \sim & \delta(\Prod_{i,j}
E(i\alpha_i\beta_j, \pi^{i+j})).
     \end{array}
     $$

To compare the left-hand side and the right-hand side in the last row, one must note that
$$
     A\partial B = \sum_{i, j} j\alpha_i\beta_jt^{i+j-1}
     ~\mbox{and}~\mbox{res}(A\partial B) = \sum_{i+j=0}j\alpha_i\beta_j.
     $$
We see that both methods are completely {\em parallel}. To be specific, the operation $\delta$ plays the role of the residue. But the second construction in the number field $K$ is, of course, much more complicated. In particular, the numbers $E(\alpha)$ and $E(\beta)$ have disappeared from our table. It is natural to compare their role in the definition of the Hilbert symbols with the so-called tame symbol in the field $\C((t))$ rather than with the residue at the point $P$ \footnote{If $f$ and $g$ belongs to $K^{*}$ then the tame symbol $(f, g)$ can be defined as $(-1)^{mn}f^{-n}g^m (t = 0)$ where the $m$ and $n$ are  valuations of the functions $f$ and $g$ in  the point $P$.}.

\section{ The General Situation in the 1950s and 1960s}

The 1950s were a period of reawakened interest in algebraic geometry in the Soviet Union (although it may not be quite accurate to speak of a "reawakening", since up to the 1950s no one in the USSR was interested in algebraic geometry)\footnote{One can mention only N. G. Chebotarev and,in particular,his book \cite{Ch}, and the papers of
I. G. Petrovskii and O. A. Oleinik on the topology of real algebraic varieties, written just after the war. See their survey in \cite{O}.}.

Nevertheless, after the war there were several people seriously interested in algebraic geometry, first among them I.R.Shafarevich, who tried to study the available literature. To show how difficult this was to do, say, in the late 1940s, there was a seminar at Moscow University in which several mathematicians, Shafarevich among them, attempted to understand the proof of the Mordell-Weil theorem, but were unable to do so.

One cause of this situation is understandable: the strict isolation from the rest of the world. For example, when mathematicians from the whole world met in the International Congress at Cambridge, Massachusetts, in 1950 for the first time after the war, there was no one there from the USSR. There was only a telegram communicating that "Soviet mathematicians, who are extremely busy with their current research, cannot participate in the Congress" \cite{Cam}. This was the very Congress at which A.Weil gave the lecture we mentioned above. In the late 1950s the situation improved somewhat, but, of course, strong restrictions remained.

The few visits from western mathematicians, among whom one must mention Erich
K\"ahler, exerted a great influence on the development of ideas during the 1950s. Because of the rarity of direct contacts, the study of the literature was very important. So far as I know, the notes of the H. Cartan seminar \cite{Car} were difficult to find in Moscow, but they were studied very thoroughly. The book of Hodge on harmonic integrals \cite{H} and the lecture notes of Siegel on automorphic functions of several complex variables \cite{Si1} were also very popular. The latter were translated into Russian in 1954 by I. I. Pyatetskii-Shapiro \cite{Si2} and in the mid-1950s Shafarevich and Pyatetskii-Shapiro conducted a seminar on this book, which, from the point of view of understanding the proofs of the theorems in the book, turned out to be much more successful. Perhaps the work of Pyatetskii-Shapiro on bounded domains and his solution of Cartan's problem on the existence of nonsymmetric bounded domains were the result of this activity. (See his reminiscences of this time in the collection \cite{G}.) In1960 and 1961 a large seminar on the theory of deformations of complex structures, which had recently been created by K. Kodaira and D. Spencer, was organized in Moscow University by E. B. Dynkin,M. M. Postnikov, andI. R. Shafarevich.

But the most important aspect for our history is the interest in the classical theory of algebraic surfaces. This is understandable from the point of view of the analogy discussed above. The constructions described in Section 1 belong to classical algebraic number theory and thus belong to dimension 1 according to our classification. Shafarevich later began to study Diophantine equations, in particular, elliptic curves of algebraic number fields, and realized the necessity of moving on to higher dimensions. After all, schemes of dimension 2 must correspond to curves over such fields. Fortunately, the concept of a scheme itself had only just arisen. This program was concisely formulated in his Stockholm lecture, to which we referred earlier. To understand arithmetic in dimension 2 one must first have a clear picture for algebraic surfaces. That is, we must have a theory of the corresponding geometric objects Ч algebraic surfaces, primarily over the field of complex numbers and then over other fields.

Such a theory already existed in the work of Italian mathematicians who studied algebraic geometry Ч G.Castelnuovo, F.Enriques, F.Severi and others. But the main definitions and proofs of the Italian geometers were not sufficiently rigorous and were sometimes simply incomprehensible. In fact, this subject was a rather isolated area of mathematics, having its own rules and laws, which were rejected by the greater part of the mathematical community. The rise of the theory of sheaves, which came out of complex-analytic geometry (the paper of J.-P.Serre \cite{S}) and the analytic methods in the papers of Kodaira and Spencer, made it possible to give new, rigorous proofs of many results of the Italian school. It suffices to compare the lecture of B. Segre at the 1954 International Congress in Amsterdam, which belongs entirely to the earlier epoch, with the lecture of Grothendieck at the 1958 Congress in Edinburgh, to get a sense of the revolution that had taken place.

In the early 1960s Shafarevich organized a seminar at Moscow University, in which the classical works of the Italian mathematicians on the theory of algebraic surfaces were studied. The main source was the book of Enriques \cite{En}. This seminar was conducted in two stages, during the 1961-62 and 1962-63 academic years. It is interesting that about the same time (more precisely, in the late 1950s) interest in the results of the Italians in the area of algebraic surfaces also arose in the USA,in the schools of O. Zariski and K. Kodaira.

The result of the two-year study was the publication of the book {\em Algebraic Surfaces} \cite{AS}, which appeared in 1965 in the Proceedings of the Steklov Mathematical Institute. This volume contained the following chapters:
\begin{enumerate}
     \item Birational transformations (A. B. Zhizhchenko)
     \item Minimal models (A. B. Zhizhchenko)
     \item Rationality criteria (A. B. Zhizhchenko)
     \item Ruled surfaces (I. R. Shafarevich)
     \item Minimal models of ruled and rational surfaces (Yu. I. Manin, A. N. Tyurin, Yu. R. Vainberg)
     \item Surfaces of general type (B. G. Moishezon)
     \item Surfaces with a pencil of elliptic curves (I. R. Shafarevich)
     \item Algebraic surfaces with  $\kappa = 0$ (B. G. Averbukh)
     \item The space of moduli of complex surfaces with $q = 0$ and $K = 0$ (G. N. Tyurina)
     \item Enriques surfaces (B. G. Averbukh).
\end{enumerate}

The book gave a complete exposition of the classification of algebraic surfaces, as it had been done by the Italians, with results proved in the language of sheaves. In some places the classical propositions were corrected or extended. Remarkably, this seminar and the book served as the main impetus for the further development of algebraic geometry in Moscow. We give here just a short list of the further research that grew out of it:
\begin{itemize}
\item rational surfaces and multidimensional varieties (with the solution of
Luroth's problem
\footnote{
Shafarevich heard the statement of this problem from Chebotarev, who had been interested in it for some time. In particular, Chebotarev had discussed the problem in his lecture at the Z\"urich Congress in 1932. The problem is to explain how a subfield of the field of rational  functions $k(x_1,\dots,x_n)$ in $n$ variables can again be a field of the same type. This is true for $n = 1$ and for $n = 2, k = \C$. The proof of this last fact was given by the Italians and used
the full power of the theory of surfaces. Manin and Iskovskikh constructed counterexamples in dimension 3. Independently, the problem was solved by P. A. Griffiths and H. C. Clemens, M. Artin and D. Mumford.}
and the classification of Fano varieties)(Yu. I. Manin and V. A. Iskovskikh)
\item the theory of vector bundles on algebraic curves and surfaces(A. N. Tyurin, F. A. Bogomolov)
\item K3 surfaces (G. N. Tyurina, I. R. Shafarevich, I. I. Pyatetskii-Shapiro,\\
V. V. Nikulin, A. N. Rudakov, V. S. Kulikov, and others)
\item elliptic sheaves and the main homogeneous spaces (I. R. Shafarevich,
O. N. Vvedenskii)
\item multidimensional birational and analytic geometry (including criteria for ampleness)(B. G. Moishezon)
\item  minimal models for arithmetical surfaces(I. R. Shafarevich).
\end{itemize}

The majority of papers in this list were written after the seminar and under its influence. The papers of Shafarevich on the theory of the principal homogeneous spaces, which preceded the seminar, are an exception. They arose out of his interest in Diophantine equations, primarily the theory of elliptic curves. As early as 1956, in a lecture given at the third All-Union Mathematical Congress, Shafarevich pointed out the analogy between embedding problems in the Galois theory of algebraic number fields and the classification problem for elliptic curves defined over such fields. What these problems had in common was their statement in the languages of Galois cohomology and the presence of local invariants connected with the completions of the base number field. (For more details, see \cite{Dem}.) It was natural to pass from these arithmetical problems to the study of elliptic curves over a field of algebraic functions, and that is what the surfaces with a pencil of elliptic curves in the preceding list amount to.

A more detailed exposition of the subsequent development of algebraic geometry in Moscow can be found in \cite{ALNT, Dem, Isk}. The general atmosphere of this era is well shown in \cite{G}.

\section{Arithmetical Surfaces}

The development of the last area in our list was of great significance for number theory. In his 1966 lectures in Bombay (now Mumbai) \cite{Sh3} Shafarevich gave a systematic development of the fundamental concepts and results from the theory of algebraic surfaces for the case of arithmetical surfaces. In these lectures, using the language of schemes, he constructed a theory of intersections, defined and studied birational transformations and minimal models \footnote{Some of these results were obtained independently by S.Lichtenbaum \cite{Lich}[34] in the USA and
P. Deligne \cite{D1} in France.}.

As an illustration, we give the simplest example of an arithmetical surface arising from the affine line over the field $\Q$. Let $X = \Spec \Z[t]$. This is a scheme of dimension 2, and it is mapped onto $B = \Spec(\Z)$.
The fibers of this mapping over points $p \in S$ are affine lines over the finite fields $\Fp$.
\par\smallskip
\includegraphics{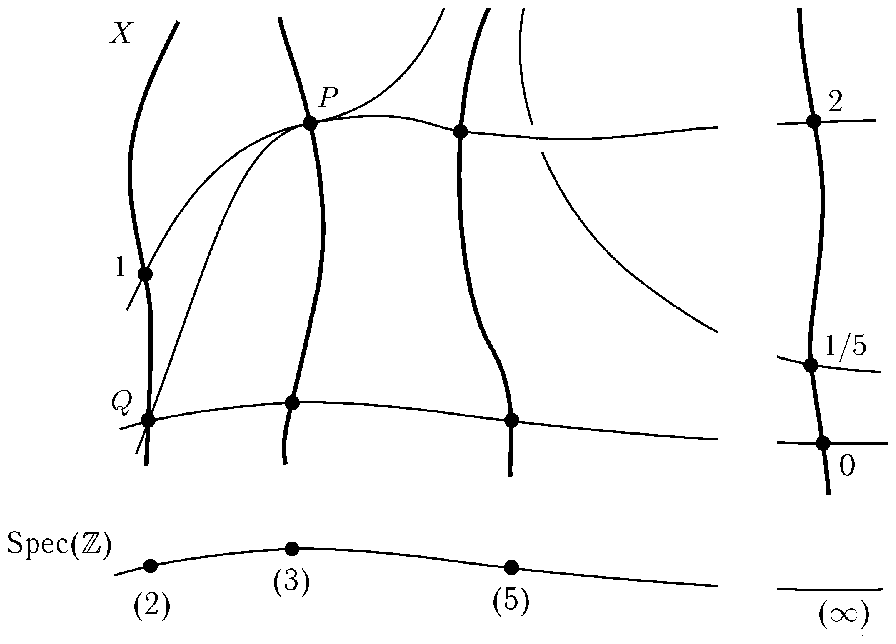}
\par\smallskip
Here we represent the points of the fibers(which are simultaneously points of the scheme $X$) with coordinates in finite fields (that is, residues mod $p$). The "surface" $X$ contains "curves" defined by equations of the form $f = const$, where $f \in \Z[t]$.
The curves $t = 0$ and $t = 2$ intersect at a point $Q$ of the fiber over $p = 2$ and have first-order tangency there, that is, they are transversal. The curves $t = 1/5$ and $t = 2$ intersect at a point $P$ of the fiber over $p = 3$ and have second-order tangency. Indeed,
$$
     2 \equiv 0~\mbox{mod}~2,~2 \not\equiv 0~\mbox{mod}~2^2
     $$
     $$
     5\cdot 2 \equiv 1~\mbox{mod}~3,~5\cdot 2 \equiv 1~\mbox{mod}~3^2,
     5\cdot 2 \not\equiv 1~\mbox{mod}~2^3.
     $$
In the last case, in a neighborhood of the fiber over $p= 3$ we have the 3-adic series expansions
$$
     2 = 2 + 0\cdot 3 + 0\cdot 3 + \dots,
     $$
     $$
     1/5 = 2 + 0\cdot 3 + 1\cdot 3^2 + \dots
     $$

The general definition of the index of intersection of two curves $C =(f = 0)$ and $D =(g = 0)$ at a point $x$ looks as follows:
\begin{equation}
\label{i}
(C\dot D)_x = \mbox{log}\#k(x) \cdot \mbox{length}\Z[t]/(f,g),
         \end{equation}
 where $\mbox{log}\#k(x)$
 is introduced in analogy with the one-dimensional case (see Section 1). Of course, this definition makes sense only if the curves $C$ and $D$ intersect in a finite set of points, that is, have no common components. To give the definition in the general case in algebraic geometry one normally uses the shift method, bringing the curves into general position. As a shift one uses the addition with the divisor of a rational function, since the index of intersection of any curve with a divisor of a function on a complete surface is zero. This last property is a generalization of a property of divisors of functions on curves (their degree is zero, see Section 1) to the case of a surface. As we have seen, for this last property to hold one must have a complete or compact curve.

Accordingly, in the two-dimensional situation one must have something like a complete surface. However, only incomplete schemes defined over an affine base (the spectrum of the ring of integers of the field of algebraic numbers) were considered in Shafarevich's lectures. It was clear from the very beginning that such an approach is only a partial analog of the situation with algebraic surfaces. At the end of the lectures the problem was posed: to find a complete analog of an algebraic surface and construct a theory of intersections for it. Let us consider this problem in more detail.

Comparing the geometric and arithmetical cases in dimension 1, we saw that the complete analog of a projective curve $X$ is the set $X=X' \bigcup \infty$ and the
structure of a scheme is present only on the subset $X'$. The point $\infty$  is adjoined
to X Уby hand,Ф so to speak. It is unknown what structure must be on the entire set $X$. It seems that the theory of schemes is unsuitable for this purpose\footnote{Recently, N. Durov from  Petersburg has suggested a generalization of the scheme theory that will cover the case also (See Nikolai Durov, {\em New Approach to Arakelov Geometry},
arXiv(math.AG): 0704.2030).}.

This is also true for higher dimensions. The complete object on the right-hand side of the table above, which corresponds to projective algebraic surfaces consists of the arithmetical surfaces introduced by S.Yu.Arakelov \cite{Ar1, Ar2} in the early 1970s.
It is not very convenient to compare them directly with algebraic surfaces. For such a comparison an algebraic surface $X$ is usually endowed with the structure of a pencil of algebraic curves parameterized by a projective nonsingular curve $B$. Thus we have a mapping
\begin{displaymath}
f:X \longrightarrow B,
\end{displaymath}
whose fibers are projective curves and almost all them are nonsingular curves of the same genus $g$. Let us now compare this mapping with the structure mapping
\begin{displaymath}
f:X' \longrightarrow Spec(\Z )
\end{displaymath}
for an arithmetical surface. Since the basis "curve" $\Spec(\Z)$ is not complete, this means that the two-dimensional scheme $X$ is also non-complete, and thus cannot be regarded as a precise analog of the algebraic surface $X$. In exactly the same way as in the case of dimension 1, we need to adjoin something.

To understand Arakelov's idea, let us return to the case of an algebraic surface $X$ with the mapping $f$ and divide the basis curve $B$ into two distinct parts $B'$ and $S$, where $B'$ is an open subset and $S$ is a finite subset. We wish
to regard $B'$ as the analog of $Spec(\Z )$
   and $f^{-1}(B')$  as the analog of $X'$. We now seek the missing part of the arithmetical surface that corresponds to the part of $X$ lying over $S$. To solve this problem we need to study this piece of the algebraic surface $X$ more attentively.

With the mapping $f$ one can connect an algebraic curve $Y$ defined over the field $K$ of rational functions of the curve $B$. In the theory of schemes this construction, which was known earlier in classical algebraic geometry, is called transition to the generic fiber and admits a simple and rigorous definition. If $b \in B$, we have the curve $Y_{(b)}$ obtained by replacing the base field $K$ by the local field $K_b$,
\begin{displaymath}
Y \otimes_{K} K_{b}
\end{displaymath}
and the two-dimensional scheme  $X_{(b)}$,
\begin{displaymath}
X \otimes_{B} Spec(O_{b}) \quad .
\end{displaymath}
obtained by replacing the base curve $B$ by an УinfinitesimalФ neighborhood $\Spec(O_b)$ of the point $b$.
Now let $b \in S$. We can then compare the field $K_{b}$ with the fields that are the completions of the field of algebraic numbers at Уinfinity.Ф In the arithmetical case we have no analogs for the schemes $X_{(b)}$, but we can define the curves $Y_{(b)}$ by the same formula as above. For the field $\Q$ this has the following appearance
\begin{displaymath}
Y_{\infty} = Y \otimes_{\Q} {\R} \subset Y \otimes_{\Q} {\C}.
\end{displaymath}
Thus we obtain Riemann surfaces over the field $\C$. Arakelov assumed that the choice of some Hermitian metric on the Riemann surfaces $Y_{\infty}$ can be regarded as replacing the nonexistent model $X_{\infty}$. Such an approach can be explained as follows. In the geometric case for the curves $Y_{(b)}$, there is a bijective correspondence
$$
  \{ \mbox{sections of the projection}~X_{(b)} \rightarrow \Spec (O_b) \}
  \longleftrightarrow Y_{(b)}(K_b)
$$
between sections of the mapping $f$ and rational points of the curve $Y_{(b)}$ (see Fig.2 below). For any two distinct sections $C$ and $D$ their index of intersection is defined (see (1)) and can be used as a metric on the curve $Y_{(b)}$. The choice of a different model $X_{(b)}$ gives another metric on $Y_{(b)}$. Thus, one can try to regard the set of models $X_(b)$ as a certain set of metrics on $Y_{(b)}$. Such an approach to the interpretation of Arakelov's theory arose much later \cite{D2}. The description of an exact correspondence between models and metrics was given only in \cite{Zh, S2}.

We now give a table that will be more precise than the general picture given above.
\vskip 0.5cm
\begin{tabular}{|c|c|}
\hline
geometric case  & arithmetical case   \\
\hline
   projective nonsingular curve     &     spectrum of the ring  $\Lambda$
                                            of integers             \\
  $B$ with finite subset
   $S \subset B$                    &  number field $K$ and embeddings
                                                                    \\
                                    &   of $K$ into $\C $                \\
\hline
 projective algebraic       &          \\
surfaces $X$ over $\Fq $      &      arithmetical surfaces       \\
 with mapping       &                              \\
 $f:X \rightarrow B$ onto  $B$        &                                \\
\hline
 surface $X'= f^{-1}(X-S)$ with    &     two-dimensional scheme   $X'$     \\
 mapping  $f\vert_{X'}$ on $X-S$          &   over  $Spec(\Lambda)$    \\
\hline
  algebraic curve  $Y_{(b)}$ with &  compact Riemann surfaces\\
  $b \in S$                         &  $Y_{(\infty)}$, corresponding  of K  \\
                                    &  to embeddings of $K$ into $\C $          \\
\hline
  schemes   $X_{(b)}$ with $b \in S$       & Hermitian metrics on surfaces\\
                                    &   $Y_{(\infty)}$                 \\
\hline
\end{tabular}
\vskip 0.5cm
Arakelov then defined such concepts as divisor, divisor of a function and differential form, linear equivalence, index of intersection, and canonical class. He proved an analog of the adjunction formula and also stated an analog of the RiemannЦRoch theorem in\cite{Ar2}.
ArakelovТs construction lay undisturbed for nearly ten years, and only in the early 1980s it did serve as the point of departure for further development in the papers of Gerd Faltings. We refer to \cite{ZP, F1, F2, S1, S3, Sz, P2}, where these later events are related. This line exerted a powerful influence on the development of number theory and also on the development of elementary-particle physics \cite{B}, demonstrating the notorious "unreasonable effectiveness of mathematics in the natural sciences."

\section{Height and ArakelovТs Theory}

In this section we shall explain the origins of Arakelov's theory, starting from the concept of height - a basic tool of the theory of Diophantine equations.

Let $X$ be a projective algebraic variety defined over a global field $K$ of dimension 1 (in other words, $K$ is either the field of algebraic numbers or the
field $K = k(B)$ of algebraic functions on some curve $B$ defined over the field of constants). And let $D$ be a divisor on $X$ (that is, an integer linear combination of subvarieties of codimension 1).

A {\em height} is a function
$$
h_{X,D}: X(K) \rightarrow \R,
$$
on the set of rational points $X(K)$ depending on the choice of the divisor $D$ on
$X$. Actually, the height is not uniquely determined by these data. We shall write $f \approx g$ for two numerical functions if $f - g$ is a bounded function. The height is defined as an equivalence class of functions relative to such an equivalence relation. (For details, see \cite{L1}).
Here is a simple but important example. Let $X \subset {\bf P}^n, K = \Q$ and let $D$  be a hyperplane section. Then the point $P \in X(K)$ is $(x_0:\dots:x_n)
     \in  {\bf P}^n(\Q)$, where $x_i$ are relatively prime integers. We have
$$
     h_{X,D} = \mbox{max}_i~log\vert x_i \vert.
$$
From this one can see that the number of points of a bounded height is finite - a very important property, which makes it possible to obtain various kinds of finiteness theorems for Diophantine equations.
More generally, for a global field $K$ with a set of valuations $\nu$ (in which we include infinite points in the number field case) and norms $\vert \cdot \vert_{\nu}$ corresponding to them, the height of a point in projective space is defined as the product
$$
     h(P) = \prod_{\nu}\mbox{max}_i~log\vert x_i \vert_{\nu}
     $$
and the product formula (see Section 1)
     $$
     \prod_{\nu} \vert x \vert_{\nu} =  1, ~x \in K^{*}
     $$
shows that this expression is well defined (but depends, for example, on the choice of the system of projective coordinates). The height has the following properties:

i) {\sc Invariance under linear equivalence}: if $D' = D + (f)$, where $(f)$ is the divisor of the function $f$, then
$$
                 h_{X, D'} \approx h_{X, D}.
$$

ii) {\sc Functoriality}: if $f: X \rightarrow Y$ is a morphism of algebraic varieties and $D$ is a divisor on $Y$, then
$$
        h_{X,f^{*}D} \approx f^{*}h_{X,D},
$$
where $f^{*}$
  is the inverse image of the divisors or functions respectively.

iii) {\sc Finiteness}: if the divisor $D$ is a hyperplane section, then for every $h \in R$ the set
$$
    \{P \in X(K) :  h_{X,D}(P) \le h \}
$$
is finite.

Using these properties A.Weil proved that the group $A(K)$ of rational points on an Abelian variety $A$ over a global field $K$ is finitely generated\footnote{Modulo the $K/k$-trace in the geometric case. See \cite{L1}.}(the Mordell-Weil theorem).

If we are in the geometric situation (according to the preceding classification), the base field $K$ has the form $k(B)$, and there exist a projective variety $Y$ and a morphism $f: Y \rightarrow B$ with a general fiber $X$. Under reasonable hypotheses on $X$ and $Y$ (irreducibility and flatness of the morphism $f$) there is a bijective correspondence
$$
       \{   \mbox{section of the mapping}   f:   Y   \rightarrow   B    \}
\longleftrightarrow X(K)
     $$
between sections $C$ and rational points $P$ on $X$. The divisor $D$ also defines a certain divisor $\bar D$ on $Y$
\par\smallskip
\includegraphics{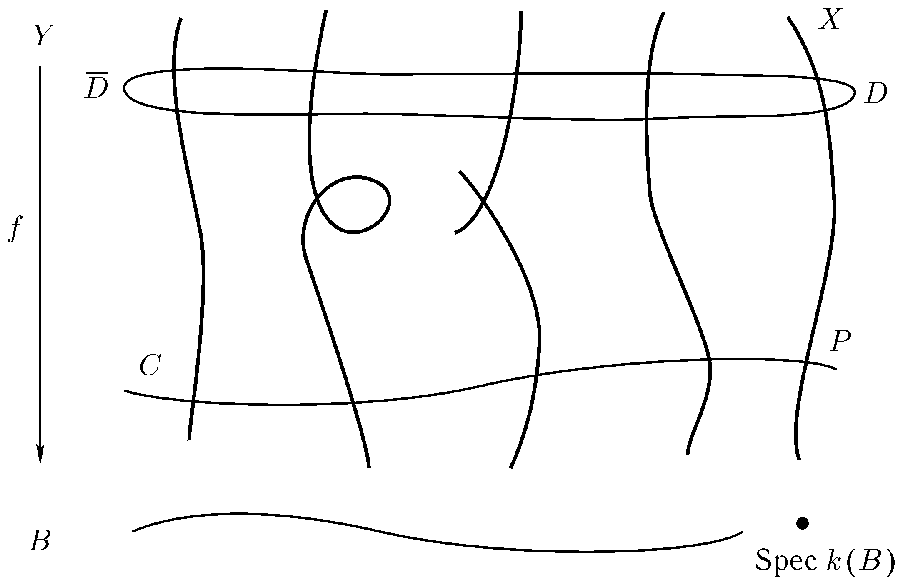}
\par\smallskip
Under these conditions we have
\par\smallskip
     iv) $h_{X,D} \approx (C.D)_Y$,
\par\smallskip
(For this equality to make sense it is, of course, necessary to have a theory of intersections on the variety $Y$, for example, to assume that $Y$ is a nonsingular variety.)

Thus, in essence, the height is the index of intersection and this circumstance can be used in different directions. If $X$ is an algebraic curve, then $Y$ is a surface, its model from the preceding section. This connection suggests that the height may serve as a starting point for a construction of a theory of intersections on arithmetical surfaces.

An obvious defect of the height is the approximativeness of its definition and its functoriality, only up to the equivalence relation indicated above. J.Tate invented a new definition of height on Abelian varieties $A$. This is a canonical function $\hat{h}_{A, D}$ on the set of rational points that behaves functorially relative to homomorphisms of Abelian varieties and is such that $\hat{h}_{A, D} \approx h_{A, D}$.

For the index of intersection in the geometric case we have the obvious expansion
$$
     C._{Y}D = \sum_{b \in B} C._{b}D
     $$
over the indices of intersection in all the fibers of the mapping $f$. A. Neron \cite{N}
found a new approach to the construction of Tate's height on Abelian varieties, which simultaneously gave a local expansion for it over points of the base $B$ (or valuations $\nu$ of the base field $K$):
$$
     \hat{h}(P)  = \sum_{\nu}h_{\nu}(P) + \sum_{\infty}h_{\infty}(P),~P \in
A(K).
     $$
We remark that in contrast to the global function $\hat{h}(P)$ the local components are not defined for all $P\in A(K)$ but only for $P \in (A - D)(K)$, becoming infinite on the divisor $D$. Thus, they do depend on the divisor $D$, not
only its linear equivalence class.

In the number field case the definition of local components is quite varied, depending on the nature of the point (nonarchimedean $\nu$ or archimedean $\infty$). For $\hat{h}_{\nu}(P)$ one uses the index of intersection on a special nonsingular model of an Abelian variety $A$ over the ring of integers of the base field $K$ (Neron's minimal model).

Analysis first enters the game at infinity. Let $A$ be an Abelian variety over the field of complex numbers $\C$, and let $D$ be a positive divisor on $A$ (that is, $D = \sum_i n_iD_i$, where $n_i \in \Z, n_i \ge 0$ and $D_i$ is an irreducible subvariety of codimension 1). Then $A$ is a complex torus $\C^n/\Gamma$, where $\Gamma$ is a discrete subgroup of rank $2n$ in $\C^n$.

On $A$ itself the divisor $D$ is not the divisor of poles of any holomorphic function, but one can find such a function on the universal covering $\C^n$ . And, what is important, this function can be constructed in a canonical manner.
The divisor $D$ is an algebraic cohomology class, that is, an element of the space $H^{1,1} \subset H^2(A, \C)$ according to the Hodge decomposition in cohomology of the variety A. On an Abelian variety the space $H^{1,1}$ consists of Hermitian matrices. If $H$ is the matrix corresponding to the divisor $D$, there exists a unique (suitably normalized) theta-function $\theta (z) =
\theta_D(z)$) on $\C^n$ having the following properties:
\par\smallskip
i) the divisor of the poles of $(\theta_D)$is $D$;
\par\smallskip
ii) $\theta(z+\gamma) = \theta(z)\mbox{exp}(\pi H(\gamma,z)  +  \pi/2
H(\gamma,\gamma))\cdot \chi(\gamma),$
\par\smallskip
where $z \in \C^n, \gamma \in \Gamma$ and $\vert \chi(\gamma)\vert = 1$
(for details, see \cite{W2}).

The local component of the height $\hat{h}_{\infty}(P)$)can now be defined as follows \cite{N}:
$$
     \hat{h}_{\infty}(z) = - \mbox{log}\vert \theta_D(z)\vert + \pi H(z,z).
     $$
It follows from property ii) of the theta-function that this function is invariant relative to $\Gamma$, that is, it is a function of the point $P \in A(\C)$. Moreover, locally on $A$, in a neighborhood of each point of the divisor $D$ we have
\begin{equation}
     \label{q}
\hat{h}_{\infty}(P)  \approx \mbox{log}\vert \mbox{holomorphic equation for}~D \vert.
     \end{equation}
One can now look at NeronТs construction from a different point of view. How can the function $\hat{h}_{\infty}$ be distinguished among all the smooth real-valued functions on $(A - D)(\C)$ satisfying ~(\ref{q})? We remark that all functions having property ~(\ref{q}) differ from one another by a function that is bounded and smooth on all of $A$.

It is not difficult to see that the condition that distinguishes $\hat{h}_{\infty}$ is the Poisson equation
\begin{equation}
     \label{h}
     \Delta \hat{h}_{\infty}  = const,~\mbox{outside}~D,~\mbox{or}~\delta_D~
     \mbox{over all}~A.
     \end{equation}
Here
     $$
     \Delta =\sum_{i,j} \frac{\partial^2}{\partial z_i\partial \bar{z}_j}
     $$
 is the Laplacian corresponding to the flat metric on $\C^n$, which, being $\Gamma$-invariant, can be lowered to $A$. $\delta_D$ is the delta-function corresponding to the subvariety $D$. This observation suggests that the definition of the local components $\hat{h}_{\infty}(P)$ can be given for any algebraic variety $X$ if one fixes some Hermitian metric on it. Then for every divisor $D$ there exists a function $\hat{h}_{\infty,D}$ satisfying conditions
(\ref{q}) and (\ref{h}) that is unique (up to a constant). Such a definition was introduced by the author in \cite{P1} and served as the point of departure for the development of Arakelov's theory.

We can now describe Arakelov's theory as follows. An arithmetical surface $X$ consists of a nonsingular two-dimensional scheme $X$ and a mapping of it
$f: X \rightarrow B$ onto the one-dimensional scheme $B = \Spec(\Lambda_K)$), where $\Lambda_K$ is the ring of integers of the field $K$ of algebraic numbers. We denote the set of infinite points of the field $K$ by $B_{\infty}$, and for each $v \in B_{\infty}$ we choose a Hermitian metric $\mu_v$  on the Riemann surface $X_v = X \otimes_v \C$.

An Arakelov divisor $\tilde{C}$ on $\tilde{X}$ a linear combination of an ordinary divisor $C$ on $X$ and the УinfiniteФ fibers $X_v$, and the latter occur with real coefficients. Let
$$
     \tilde{C} = C + \sum_{v \in B_{\infty}} a_{v}X_v,~
     \tilde{D} = D + \sum_{v \in B_{\infty}} b_{v}X_v,~
     $$
be two Arakelov divisors.  Assume that $C$ and $D$ intersect in a finite set of
points. Then
$$
     \tilde{C}\cdot \tilde{D} = C\cdot_XD   + \sum_{v \in B_{\infty}}
    (C\cdot D)_{v}  + \sum_{v \in B_{\infty}} a_{v} \mbox{deg}D/B +
 \sum_{v \in B_{\infty}} b_{v}\mbox{deg}C/B,
     $$
where $C\cdot_XD$ is the index of intersection on the scheme $X$, and the archimedean indices $(C\cdot D)_{v}$ are defined using the GreenТs functions $G(P,Q) (= G_v(P,Q))$ constructed with respect to the metric $\mu_v$.

We recall that a Green function on a Riemann surface $X = X_v$ is determined uniquely by the following conditions:
\begin{enumerate}
\item $G$  is a smooth real-valued positive function on
$(X \times X) - (\mbox{diagonal})$.
\item If $z$ is a local holomorphic coordinate near the point $P_0$  on $X$, then
near $(P_0, P_0)$    the function $G(P, Q)$   has the form
$$\vert z(P) - z(Q) \vert
\cdot(\mbox{smooth nonvanishing function}.
$$
\item  $$\Delta_Q\mbox{log}G(P, Q) = d\mu/dz\wedge d\bar{z},$$
where
$$
\Delta_Q = (1/\pi i)(\partial^2/\partial z \partial \bar{z})
$$
is the Laplacian and $d\mu$ is the volume element that arises from the metric $\mu_v$.
\end{enumerate}
\par\smallskip
Let us set
$$
(C\cdot D)_v = \sum_{P, Q} n_Pm_Q \mbox{log}G_v(P, Q),
$$
if $C = \sum n_PP$~ and~ $D = \sum m_QQ$ are the expansions of the divisors into (finite) sums of points on $X_v(\C)$.

If $F$ is a rational function on $X$ we define its Arakelov divisor as
$$
(\tilde{F})  = (F)_X + \sum_vX_v,~a_v = -\int\mbox{log}\vert F\vert
d\mu_v.
$$
Here $(F)_X$ is the usual divisor of the function $F$ in the scheme $X$.
We can now define linear equivalence $\approx$ of divisors on $\tilde{X}$:
$$
\tilde{C} \approx \tilde{D}~\mbox{if}~\tilde{C} - \tilde{D} = (\tilde{F}).
$$
The main property of the index of intersection is its invariance relative to linear equivalence
$$
\tilde{C}\cdot \tilde{D} = \tilde{C}\cdot (\tilde{D} + (\tilde{F})).
$$
This makes it possible to define the index of intersection for any two Arakelov divisors by the usual method of algebraic geometry.

Among the classes of divisors relative to linear equivalence there is,just as in ordinary geometry, a canonical class. The divisors in that class are constructed from a rational differential form $\omega$ of degree 1 on $X$ (more precisely, it is a section of the relative cotangent bundle of $X$ over $B$). We set
$$
(\tilde{\omega}) = (\omega)_X + \sum_v a_vX_V,~a_v = -\int_{X_v}\mbox{log}
\vert \frac{\omega \wedge \bar{\omega}}{d\mu_v}\vert d\mu_v,
$$
where $(\omega)_X$ is the divisor of the form $\omega$ in the scheme $X$.

The adjunction formula for the divisor $\tilde{C} = C$, which represents a section of $C$ on an arithmetical surface(that is, a divisor having degree 1 over the base $B$), has the form
$$
  \tilde{C}\cdot (\tilde{\omega}) + \tilde{C}\cdot \tilde{C} = 0.
%  \mbox{log}\vert D_{B/\Z} \vert,
$$
In ordinary algebraic geometry the adjunction formula for a curve $C$ on a surface $X$ is the following
$$
C\cdot (\Omega_X) + C\cdot C = 2g - 2,
$$
where $\Omega_X$ is the class of differential forms of degree 2 on $X$ and $g$ is the genus of the (nonsingular) curve $C$. If the surface $X$ is fibered over the curve $B$ and $C$ is a section of this fiber bundle, then $2g - 2 = C\cdot f^{*}(\Omega_B)$ (here $\Omega_B$ is the canonical class of the curve $B$), and the adjunction formula has the form
$$
C\cdot ((\Omega_X) - f^{*}(\Omega_B)) + C\cdot C = C\cdot((\Omega_{X/B})) + C\cdot C = 0,
$$
where $(\Omega_{X/B})$ is the class of divisors corresponding to the relative cotangent bundle of the surface $X$ over $B$.

It is this equality that carries over to the arithmetical surfaces of Arakelov. For a more detailed exposition of ArakelovТs theory and its subsequent development for varieties of any dimension see \cite{F1, F2, L2, S1}.
In our brief exposition we have examined only one stem on the enormous tree of the analogy between numbers and functions. Some idea of the tree as a whole can be gained from the following list:
\begin{itemize}
\item class-field theory (a parallel description of Abelian extensions of number and function fields);
\item the zeta-and $L$-functions of schemes of dimension 1 (the problem of meromorphic continuation and the proof of the functional equation);
\item the theory of height of rational points in Diophantine geometry;
\item the Arakelov theory of arithmetical varieties;
\item the classification of semi-simple algebraic groups over local fields;
\item the theory of BruhatЦTits buildings and symmetric spaces;
\item arithmetical subgroups of algebraic groups, in particular the theory of reduction;
\item the Langlands program of describing representations of Galois groups of local and global fields;
\item the analogy between explicit formulas in number theory and the Lefschetz trace formula(A.Weil, C.Deninger, A.Connes).
\end{itemize}
This list is surely incomplete\footnote{We mentioned above the influence of number theory on physics. Then the p-adic analogs of quantum mechanics and certain models of field theory were constructed. (See, for example, the survey \cite{VVZ}).}
 and the whole story is still far from over.

\end{document}